\documentclass[12pt]{amsart}
\usepackage{amsfonts}
\usepackage{amssymb}
\theoremstyle{definition}

\theoremstyle{remark}

\begin{document}

\centerline{\large\bf $AH_3$-MANIFOLDS \ OF \ CONSTANT}

\vspace{0.1in}
\centerline{\large\bf ANTIHOLOMORPHIC \ SECTIONAL \ CURVATURE
\footnote{\it PLISKA Studia mathematica bulgarica. Vol. 9, 1987, p. 52-57.}}

\vspace{0.2in}
\centerline{\large OGNIAN T. KASSABOV}

\vspace{0.4in}
{\sl The purpose of this paper is to prove that an $AH_3$-manifold of constant 
antiholomorphic sectional curvature is a real space form or a complex space form.}

\vspace{0.1in}
{\bf 1. Introduction.} Let $M$ be a $2m$-dimensional  almost Hermitian manifolds 
with metric tensor $g$ and almost complex structure $J$. The Riemannian 
connection and the curvature tensor are denoted by $\nabla$ and $R$, respectively. 

If $\nabla J=0$, or $(\nabla_XJ)X=0$ or 
$$
	g((\nabla_XJ)Y,Z)+g((\nabla_YJ)Z,X)+g((\nabla_ZJ)X,Y)=0 \,, 
$$
then $M$ is said to be a K\"ahler, or nearly
K\"ahler, or almost K\"ahler manifold, respectively. The
corresponding classes of manifolds are denoted by $K$, $NK$, $AK$. The general
class of almost Hermitian manifold is denoted by $AH$. If $L$ is a class 
of almost Hermitian manifolds, its subclass of $L_i$-manifolds is defined by
the identity $i)$, where

1) $R(X,Y,Z,U)=R(X,Y,JZ,JU)$;

2) $R(X,Y,Z,U)=R(X,Y,JZ,JU)+R(X,JY,Z,JU)+R(JX,Y,Z,JU)$;

3) $R(X,Y,Z,U)=R(JX,JY,JZ,JU)$.

\noindent
It is well known, that 
$$
	K=K_1\subset NK=NK_2\,, \qquad K\subset AK_2\,,
$$
$$
	K=NK\cap AK\,,  \qquad AH_1\subset AH_2 \subset AH_3 \,,
$$
see e.g. [4].

A plane $\alpha$ in  $T_p(M)$ is said to be holomorphic 
(resp. antiholo\-morphic) if $\alpha = J\alpha$ (resp. 
$\alpha \perp J\alpha$). The manifold $M$ is said to be of pointwise constant
holomorphic (respectively, antiholomorphic) sectional curvature $\nu$, if for 
each point $p\in M$ the curvature of an arbitrary holomorphic (resp. 
antiholomorphic) plane $\alpha$ in $T_p(M)$ doesn't depend on $\alpha$:
$K(\alpha)=\nu(p)$.

For K\"ahler manifolds the requirements for constant holomorphic and constant
antiholo\-morphic sectional curvature are equivalent [2]. In [3] it is proved a
classification theorem for nearly K\"ahler manifolds of constant holomorphic
sectional curvature.

If $M$ is a $2m$-dimensional $AH_3$-manifold of pointwise constant antiholomorphic
sectional curvature $\nu$, and if $m>2$, then $\nu$ is a global constant [5].
In [1] it is proved a classification theorem for nearly K\"ahler manifolds of constant
antiholomorphic sectional curvature and a corresponding result for $AK_3$-manifolds
is obtained in [6].

In section 3 we shall prove the following theorem:

{\bf Theorem.} {\it Let $M$ be a $2m$-dimensional $AH_3$-manifold, $m>2$. If $M$ is
of pointwise constant antiholomorphic sectional curvature, then $M$ is a real space
form or a complex space form.}

Here a real space form means a Riemannian manifold of constant sectional curvature 
and a complex space form means a K\"ahler manifold of constant holomorphic 
sectional curvature.

\vspace{0.2in}
{\bf 2. Basic formulas.} If $M$ is an $AH_3$-manifold, its Ricci tensor S satisfies
$$
	S(X,Y)=S(Y,X)=S(JX,JY)\,.
$$
If moreover $M$ has pointwise constant antiholomorphic sectional curvature $\nu$, its 
curvature tensor has the form
$$
	R=\frac 16 \psi(S) +\nu\pi_1-\frac{2m-1}3\nu\pi_2 \,,  \leqno (2.1)
$$
where
$$
	\begin{array}{r}
		\psi(Q)(x,y,z,u)=g(x,Ju)Q(y,Jz)-g(x,Jz)Q(y,Ju)-2g(x,Jy)Q(z,Ju)   \\
		+g(y,Jz)Q(x,Ju)-g(y,Ju)Q(x,Jz)-2g(z,Ju)Q(x,Jy)
	\end{array}
$$
for an arbitrary tensor $Q$ of type (0,2) and 
$$
	\pi_1(x,y,z,u)=g(x,u)g(y,z)-g(x,z)g(y,u)  \,,
$$
$$
	\pi_2=\frac12\psi(g)  \,,
$$
see [1]. According to (2.1), $M$ is an $AH_2$-manifold.

On the other hand, it is known, that if $M$ is an $AK_2$-manifold,
$$
	R(x,y,z,u)-R(x,y,Jz,Ju)=\frac12 
	        g((\nabla_xJ)y- (\nabla_yJ)x,(\nabla_zJ)u-(\nabla_uJ)z) \,,   \leqno (2.2)
$$
holds good [4].

We shall use also the second Bianchi identity
$$
	(\nabla_xR)(y,z,u,v)+(\nabla_{y}R)(z,x,u,v)+(\nabla_zR)(x,y,u,v)=0 \,.   \leqno (2.3)
$$

\vspace{0.2in}
{\bf 3. Proof of the theorem.}

{\bf Lemma.} {\it The conditions of the theorem imply that $M$ is an Einsteinian manifold.}

{\bf Proof of Lemma.} Let $p$ be an arbitrary point of $M$ and let $x,\,y \in T_p(M)$.
According to the second Bianchi identity,
$$
	(\nabla_xR)(Jx,y,y,Jx)+(\nabla_{Jx}R)(y,x,y,Jx)+(\nabla_yR)(x,Jx,y,Jx)=0 \,.   \leqno (3.1)
$$
Let $ \{ e_i,\,Je_i; \, i=1,\hdots,m \}$ be an orthonormal basis of $T_p(M)$
such that $Se_i=\lambda_ie_i$, $i=1,\hdots,m$. Putting in
(3.1) $x=e_i$, $y=e_j$ or $x=e_k$, $y=e_i+e_j$ for $i \ne j \ne k \ne i$ and
using (2.1), we obtain
$$
	(\nabla_{e_j}S)(e_i,e_j)+\{\lambda_i+\lambda_j-2(2m-1)\nu\}g(Je_i,(\nabla_{e_j}J)e_j)=0\,;  \leqno (3.2)
$$
$$
	\begin{array}{c}\vspace{0.2cm}
		 (\nabla_{e_i}S)(e_j,e_k)+\{\lambda_i+\lambda_k-2(2m-1)\nu\}g(Je_k,(\nabla_{e_j}J)e_i) \\
		+ (\nabla_{e_j}S)(e_i,e_k)+\{\lambda_j+\lambda_k-2(2m-1)\nu\}g(Je_k,(\nabla_{e_i}J)e_j)=0  \,,
	\end{array}       \leqno (3.3)
$$
respectively. Analogously from
$$
	(\nabla_{e_i}R)(Je_j,e_j,e_j,Je_k)+(\nabla_{Je_j}R)(e_j,e_i,e_j,Je_k)+(\nabla_{e_j}R)(e_i,Je_j,e_j,Je_k)=0 
$$
we find
$$
	\begin{array}{c}\vspace{0.2cm}
		\ \ 3(\nabla_{e_i}S)(e_j,e_k)+6\{\lambda_j-(2m-1)\nu\}g((\nabla_{e_i}J)e_j,Je_k) \\
		- (\nabla_{e_j}S)(e_i,e_k)-\{\lambda_i+\lambda_j-2(2m-1)\nu\}g((\nabla_{e_j}J)e_i,Je_k)=0  
	\end{array}       \leqno (3.4)
$$
and hence
$$
	\begin{array}{c}\vspace{0.2cm}
		\ \ 8(\nabla_{e_i}S)(e_j,e_k)+\{17\lambda_j-\lambda_i-16(2m-1)\nu\}g((\nabla_{e_i}J)e_j,Je_k) \\
		+3 (\lambda_i-\lambda_j)g((\nabla_{e_j}J)e_i,Je_k)=0  \,.  
	\end{array}       \leqno (3.5)
$$
In (3.5) we change $j$ and $k$ and we add the result with (3.5)
$$
	\begin{array}{c}\vspace{0.2cm}
		\ \ 16(\nabla_{e_i}S)(e_j,e_k)+17(\lambda_j-\lambda_k)g((\nabla_{e_i}J)e_j,Je_k) \\
		+3 (\lambda_i-\lambda_j)g((\nabla_{e_j}J)e_i,Je_k)+3 (\lambda_i-\lambda_k)g((\nabla_{e_k}J)e_i,Je_j)=0  \,.  
	\end{array}       \leqno (3.6)
$$
On the other hand, (3.3) and (3.4) imply
$$
	\{3\lambda_j-\lambda_i-2\lambda_k\}g((\nabla_{e_i}J)e_j,Je_k)
	  + \{3\lambda_i-\lambda_j-2\lambda_k\}g((\nabla_{e_j}J)e_i,Je_k)=0 \,.  \leqno (3.7)
$$
Hence it is not difficult to find
$$
	3(\lambda_j-\lambda_k)g((\nabla_{e_i}J)e_j,Je_k) +
	 (\lambda_i-\lambda_j)g((\nabla_{e_j}J)e_i,Je_k) +
	 (\lambda_i-\lambda_k)g((\nabla_{e_k}J)e_i,Je_j)=0
$$
and by using (3.6) this implies
$$
	2(\nabla_{e_i}S)(e_j,e_k)=(\lambda_k-\lambda_j)g((\nabla_{e_i}J)e_j,Je_k) \,.  \leqno (3.8)
$$

Let us first assume that $g((\nabla_{e_i}J)e_j,Je_k) \ne 0$. Using three times (3.7), we obtain
$$
	(3\lambda_i-\lambda_k-2\lambda_j)(3\lambda_j-\lambda_i-2\lambda_k)(3\lambda_k-\lambda_j-2\lambda_i)
$$
$$
	\ \ \ -(3\lambda_i-\lambda_j-2\lambda_k)(3\lambda_j-\lambda_k-2\lambda_i)(3\lambda_k-\lambda_i-2\lambda_j)=0
$$
or equivalently
$$
	(\lambda_i-\lambda_j)(\lambda_j-\lambda_k)(\lambda_k-\lambda_i)=0 \,.
$$
Hence it follows $\lambda_i=\lambda_j=\lambda_k$. Indeed we have to consider two cases:

C a s e \, 1. $\lambda_i=\lambda_j$. In (3.7) we made a cyclic change of $i,\,j,\,k$ and
we use $\lambda_i=\lambda_j$:
$$
	\begin{array}{c}\vspace{0.2cm}
		(\lambda_i-\lambda_k)\{3g((\nabla_{e_j}J)e_k,Je_i)+g((\nabla_{e_k}J)e_i,Je_j)\}=0 \,, \\
		(\lambda_i-\lambda_k)\{g((\nabla_{e_k}J)e_i,Je_j)+3g((\nabla_{e_i}J)e_j,Je_k)\}=0  \,.  
	\end{array}       \leqno (3.9)
$$
If $g((\nabla_{e_k}J)e_i,Je_j) = 0$ the last equation implies $\lambda_i=\lambda_k$,
i.e. $\lambda_i=\lambda_j=\lambda_k$. So we assume $g((\nabla_{e_k}J)e_i,Je_j) \ne 0$.
In (3.5) we change $i$ and $k$ and we use $\lambda_i=\lambda_j$ and (3.8):
$$
	\{ 17\lambda_i-\lambda_k-16(2m-1)\nu \}g((\nabla_{e_k}J)e_j,Je_i) +
	      3(\lambda_k-\lambda_i)g((\nabla_{e_j}J)e_k,Je_i) =0  \,.
$$
Hence, using (3.9), we obtain $\lambda_i=(2m-1)\nu$. On the other hand, (3.5)
and (3.8) result
$$
	3\lambda_i+\lambda_k-4(2m-1)\nu = 0
$$
and so we find $\lambda_k=(2m-1)\nu$, i.e. $\lambda_i=\lambda_j=\lambda_k$.

C a s e \, 2. $\lambda_j=\lambda_k$. From (3.7)we obtain
$$
	(\lambda_i-\lambda_j)\{g((\nabla_{e_i}J)e_j,Je_k)-3g((\nabla_{e_j}J)e_i,Je_k)\}=0 \,.
$$ 
If $ g((\nabla_{e_j}J)e_i,Je_k) = 0$ this implies $\lambda_i=\lambda_j$, so
$\lambda_i=\lambda_j=\lambda_k$. But $ g((\nabla_{e_j}J)e_i,Je_k) \ne 0$ is
the Case 1. 

So we have $\lambda_i=\lambda_j=\lambda_k$ and using (3.5) and (3.8), we find
$\lambda_i=(2m-1)\nu$. If $m=3$ $M$ is Einsteinian in $p$. Let $m>3$. For
$s \ne i,\,j,\,k$ we have
$$
	(\nabla_{e_i}R)(e_s,Je_s,e_j,Je_k)+(\nabla_{e_s}R)(Je_s,e_i,e_j,Je_k)+(\nabla_{Je_s}R)(e_i,e_s,e_j,Je_k)=0 \,. 
$$
Because of (2.1) this implies
$$
	(\nabla_{e_i}S)(e_j,e_k)+\{\lambda_j+\lambda_s-2(2m-1)\nu\}g((\nabla_{e_i}J)e_j,Je_k)=0\,.
$$
Hence, using $\lambda_j=\lambda_k=(2m-1)\nu$ and (3.8), we derive $\lambda_s=(2m-1)\nu$.
Consequently $M$ is Einsteinian in $p$.

Now we assume that 
$$
	g((\nabla_xJ)y,z)=0
$$
whenever $x,\,y,\,z$ are choosen among the basic vectors $e_i,\,Je_i;\,i=1,\hdots,m$ and
$x \ne y,z,Jy,Jz$. In (2.3) we put $x=Je_i$, $y=v=e_j$, $z=-Ju=e_k$ for $i\ne j\ne k\ne i$.
Using (2.1), we obtain
$$
	(\nabla_{e_i}S)(e_i,e_j)+\{\lambda_j+\lambda_k-2(2m-1)\nu\}g(Je_i,(\nabla_{e_j}J)e_j,)=0\,.
$$
From this equality and (3.2) it follows that if $g(Je_i,(\nabla_{e_j}J)e_j) \ne 0$ for some
$i,\,j$, then $\lambda_s=\lambda_k$ for $s,\,k \ne j$. Consequently if $(\nabla_{e_s}J)e_s \ne 0$
for any $s \ne j$ then $M$ is Einsteinian in $p$. 

Let us assume that $M$ is not Einsteinian in $p$. Then $M$ is not Einsteinian in a 
neighbourhoohd $U$ of $p$. We shall prove that $M$ is an $AK_2$-manifold in $U$. Let 
$q \in U$. If $M$ is a K\"ahler manifold in $q$, $M$ is an $AK_2$-manifold in $U$.
Let $M$ is not K\"ahler in $q$. Let $\{f_i,\,Jf_i,\,i=1,\hdots,m\}$ be an orthonormal
basis of $T_p(M)$, such that $Sf_i=\mu_if_i$, $i=1,\hdots,m$. Since $M$ is non K\"ahler
and non Einsteinian in $q$ we may assume that  $(\nabla_{f_1}J)f_1 \ne 0$, 
$\mu_2=\hdots =\mu_m=\mu$ and
$$
	(\nabla_xJ)y=0\,, \qquad g((\nabla_{f_1}J)x,y)=0    \leqno (3.10)
$$
whenever $x,\,y$ are choosen among $f_i,\,Jf_i$ for $i>1$. Analogously to (3.2)
$$
	(\nabla_{f_j}S)(f_i,f_j)+\{\mu_i+\mu_j-2(2m-1)\nu\}g(Jf_i,(\nabla_{f_j}J)f_j)=0    \leqno (3.2')
$$
holds good and according to (3.10) this implies
$$
	(\nabla_{f_j}S)(f_i,f_j)=(\nabla_{Jf_j}S)(f_i,Jf_j)=0  \quad {\rm for} \ \ j>1,\ j\ne i \,.   \leqno (3.11)
$$
In (2.3) we put $x=f_i$, $y=-Jv=f_j,\,z=-Ju=f_1$ for $i\ne j\ne1 \ne i$ and using 
(2.1), (3.10) and (3.11) we obtain		
$$
	\begin{array}{c}\vspace{0.2cm}
		(\nabla_{f_i}S)(f_j,f_j)+(\nabla_{f_i}S)(f_1,f_1)-(\nabla_{f_1}S)(f_i,f_1)  \\
			+2\{\mu-(2m-1)\nu\}g(Jf_i,(\nabla_{f_1}J)f_1)=0 \,.
	\end{array}      \leqno (3.12)
$$
Now let $k\ne i$. From
$$
	(\nabla_{f_i}R)(f_k,Jf_k,Jf_k,f_k)+(\nabla_{f_k}R)(Jf_k,f_i,Jf_k,f_k)+(\nabla_{Jf_k}R)(f_i,f_k,Jf_k,f_k)=0  
$$
it follows
$$
	\begin{array}{c}\vspace{0.2cm}
		2(\nabla_{f_i}S)(f_k,f_k)-(\nabla_{f_k}S)(f_i,f_k)+ \{\mu_i+\mu_k-2(2m-1)\nu\}g(Jf_i,(\nabla_{f_k}J)f_k) \\
			-(\nabla_{Jf_k}S)(f_i,Jf_k)+ \{\mu_i+\mu_k-2(2m-1)\nu\}g(Jf_i,(\nabla_{Jf_k}J)Jf_k)=0 \,.
	\end{array}      
$$
Hence using $(3.2')$ we derive
$$
	(\nabla_{f_i}S)(f_k,f_k)=(\nabla_{f_k}S)(f_i,f_k)+(\nabla_{Jf_k}S)(f_i,Jf_k)  \,.  \leqno (3.13)
$$
Now (3.11) and (3.13) imply
$$
	(\nabla_{f_i}S)(f_j,f_j)=0  \quad {\rm for} \ \ i,\,j>1,\,i\ne j \,.
$$
Then (3.12) takes the form 
$$
		(\nabla_{f_i}S)(f_1,f_1)-(\nabla_{f_1}S)(f_i,f_1)+2 \{\mu - (2m-1)\nu\}g(Jf_i,(\nabla_{f_1}J)f_1)=0 
$$
and using (3.13), we obtain
$$
		(\nabla_{Jf_1}S)(f_i,Jf_1)+2 \{\mu - (2m-1)\nu\}g(Jf_i,(\nabla_{f_1}J)f_1)=0 \\
$$
which implies
$$
	\begin{array}{c}\vspace{0.2cm}
	(\nabla_{f_1}S)(f_i,f_1)+(\nabla_{Jf_1}S)(f_i,Jf_1)   \\
	+2 \{\mu - (2m-1)\nu\}g(Jf_i,(\nabla_{f_1}J)f_1+(\nabla_{Jf_1}J)Jf_1)=0  \,.
	\end{array}	          \leqno (3.14)
$$
Since $M$ is not Einsteinian in $q$ the first equation of (3.2) and (3.14) result
$$
	(\nabla_{f_1}J)f_1+(\nabla_{Jf_1}J)Jf_1=0  \,.  \leqno (3.15)
$$
From (3.10) and (3.15) it follows easily that $M$ is an almost K\"ahler maniflod 
in $q$. Consequently it is an almost K\"ahler manifold in $U$ and hence an $AK_2$-manifold
in $U$. If $M$ is a K\"ahler manifold in $U$ it is of constant holomorphic sectional
curvature [2] and hence Einsteinian in $U$ which contradicts our assumption. Let $M$
is non K\"ahler in $q$ (we shall use the above notations for the basis of $T_q(M)$) and let
$$
	(\nabla_{f_1}J)f_i=\alpha_if_1+\beta_iJf_1  \quad {\rm for} \ \ i>1\,.
$$
In (2.2) we put $x=u=f_i$, $y=z=f_1$:
$$
	\nu-\frac16(\mu+\mu_1)+\frac{2m-1}3\nu=-\frac12(\alpha_i^2+\beta_i^2) 
$$
for $i>1$ which implies
$$
	\alpha_i^2+\beta_i^2 =\alpha_j^2+\beta_j^2   	\quad {\rm for} \ \ i,\,j>1\,. \leqno(3.16)
$$
Now we put in (2.2) ($x=f_i,\,y=z=f_1,\,u=f_1$), ($x=f_i,\,y=z=f_j,\,u=Jf_j$) 
respectively and we obtain
$$
	\begin{array}{c}\vspace{0.2cm}
		\alpha_i\alpha_j+\beta_i\beta_j=0 \,,  \\
		\alpha_i\beta_j-\alpha_j\beta_i=0 \,,
	\end{array}   \leqno (3.17)
$$ 
respectively. But (3.16) and (3.17) imply $\alpha_i=\beta_i=0$ for $i>1$ which
is a contradiction. This proves the Lemma.

Now we prove the Theorem. Since $M$ is Einsteinian (2.1) takes the form
$$
	R=\nu\pi_1+\lambda\pi_2
$$
with a constant $\lambda$. Consequently $M$ is a real space form or a complex space form [7].

\vspace{0.6in}
\centerline{\large R E F E R E N C E S}

\vspace{0.2in}
\noindent

\noindent
1. G. G\,a\,n\,c\,h\,e\,v,  O. K\,a\,s\,s\,a\,b\,o\,v. Nearly K\"ahler manifolds of 
constant antiholomorphic 

\ sectional curvature. {\it C. R. Acad. bulg. Sci.}, {\bf 35}, 1982, 145-147.

\noindent
2. B.-Y. C\,h\,e\,n, K. O\,g\,i\,u\,e. Some characterizations of complex space forms.
{\it Duke Math. 

\ J.,} {\bf 40}, 1973, 797-799.

\noindent
3. A. G\,r\,a\,y.  Classification des vari\'et\'es approximativement k\"ahleriennes de courbure sec-

\ tionnelle holomorphe constante. {\it C. R. Acad. Sci. Paris, S\'er. A}, {\bf 279}, 1974, 797-800.

\noindent
4. A. G\,r\,a\,y. Curvature identities for Hermitian and almost hermitian manifolds. 
{\it T\^ohoku 

\ Math. J.,} {\bf 28}, 1976, 601-612.

\noindent
5. O. K\,a\,s\,s\,a\,b\,o\,v. Sur le th\'eor\`eme de F. Schur pour une vari\'et\'e presque 
hermitienne. 

\ {\it C. R. Acad. bulg. Sci.,} {\bf 35}, 1982, 905-908.

\noindent
6. O. K\,a\,s\,s\,a\,b\,o\,v. Almost K\"ahler manifolds of constant antiholomorphic sectional curva-

\ ture. {\it Serdica}, {\bf 9}, 1983, 373-376.

\noindent
7. F. T\,r\,i\,c\,e\,r\,r\,i, L. V\,a\,n\,h\,e\,c\,k\,e. Curvature tensors on almost Hermitian manifolds.

\ {\it Trans. Amer. Math. Soc.}, {\bf 267}, 1981, 365-398.

\vspace {0.5cm}
\noindent
{\it Center for mathematics and mechanics \ \ \ \ \ \ \ \ \ \ \ \ \ \ \ \ \ \ \ \ \ \ \ \ \ \ \ \ \ \ \ \ \ \
Received 15.VIII.1983

\noindent
1090 Sofia   \ \ \ \ \ \ \ \ \ \ \ \ \ \ \ \ \  P. O. Box 373}

\end{document}